\documentclass[main.tex]{subfiles}

\newcommand{\kat}[1]{{\textcolor{black}{#1}}}
\usepackage[normalem]{ulem}
\newcommand{\katso}[1]{{\textcolor{red}{\sout{}}}}

\usepackage{hyperref}
\hypersetup{
    colorlinks=true,
    linkcolor=black,
    }

\title{\kat{Nowhere-zero 8-flows in cyclically 5-edge-connected, flow-admissible signed graphs}}
\author{Matt DeVos \qquad Kathryn Nurse \quad Robert \v S\'amal}
\author{
   Matt DeVos\thanks{
     Email: {\tt mdevos@sfu.ca}. 
     Supported by an NSERC Discovery Grant (Canada)}
 \and
   Kathryn Nurse\thanks{
  Email: {\tt knurse@sfu.ca} }
  \and
    Robert \v S\'amal\thanks{
    Email: {\tt samal@iuuk.mff.cuni.cz}.     
    Supported by grant 22-17398S of the Czech Science Foundation. 
    }
}

\begin{document}
\maketitle

\begin{abstract}
    \kat{In 1983, Bouchet proved that every bidirected graph with a nowhere-zero integer-flow has a nowhere-zero 216-flow, and conjectured that 216 could be replaced with 6. This paper shows that for cyclically 5-edge-connected bidirected graphs that number can be replaced with 8.}
\end{abstract}




\section{Introduction}

A \emph{signed graph} consists of a graph $G$ together with a signature $\sigma : E(G) \rightarrow \{-1,1\}$.  The \emph{sign} of a cycle $C$ is defined to be $\sigma(C) = \prod_{e \in E(C)} \sigma(e)$, and we call $C$ \emph{negative} \emph{(positive)} if $\sigma(C) = -1$ $(\sigma(C) = 1)$.  In the setting of signed graphs we are standardly interested only in the signs on cycles; accordingly, we consider two signatures to be \emph{equivalent} if they have the same sign on every cycle.  If $S$ is an edge-cut of $G$, then we may modify $\sigma$ by switching its values on the edges in $S$ to obtain an equivalent signature, and we call this operation a \emph{switch}.  It is a useful exercise to show that whenever two signatures are equivalent, they are related by a switch.  We call the signed graph \emph{balanced} if there are no negative cycles.  

Following Bouchet, we will treat every edge as composed of two half-edges.  We define $H(G)$ to be the set of half-edges of $G$, so for an edge $e$ with ends $v,v'$ we view $e$ as consisting of two half edges, say $h,h'$, where $h$ $(h')$ is incident with $v$ $(v')$.  An \emph{orientation} of a signed graph is a function $\tau : H(G) \rightarrow \{-1,1\}$ with the property that every edge $e$ composed of half-edges $h,h'$ satisfies $\tau(h) \tau(h') = - \sigma(e)$.  If $h$ is incident with the vertex $v$ then we view $h$ as directed \emph{away from} \emph{(toward)} $v$ if $\tau(h) = 1$ $(\tau(h) = -1)$.  Note that orienting a positive edge directs one half-edge toward its end and the other away (akin to an orientation of an ordinary (not signed) graph).  On the other hand, when a negative edge is oriented, both half-edges will be directed toward or away from their ends.  A signed graph that is additionally equipped with an orientation is called a \emph{bidirected} graph.  

Bouchet introduced flows on signed graphs as a non-orientable generalization of Tutte's flow-colouring duality for planar graphs.  Following his setup, we let $G$ denote a bidirected graph, $A$ an abelian group, and we let $\phi : E(G) \rightarrow A$.  We define the \emph{boundary} of $\phi$ to be the function $\partial \phi : V(G) \rightarrow A$ given by the following rule.  Here we use the notation that $H_v$ denotes the set of half-edges incident to a vertex $v$, and for a half-edge $h$ we let $e_h$ denote the edge containing it.  

\[ \partial \phi (v) = \sum_{h \in H_v} \tau(h) \phi(e_h) \]

We call the function $\phi$ an $A$-\emph{flow} if $\partial \phi$ is identically 0, and we call it a $k$-flow for a positive integer $k$ if $A = {\mathbb Z}$ and $\phi(E(G)) \subseteq \{ - (k-1), \ldots, (k-1) \}$.  We say that the function $\phi$ is \emph{nowhere-zero} if $0 \not\in \phi(E(G))$.  Note that whenever $\phi$ is a flow and $e$ is an edge made of the half-edges $h,h'$,  we may modify our orientation by switching the values of $\tau(h)$ and $\tau(h')$ and modify $\phi$ by replacing $\phi(e)$ with its negation to obtain a new flow.  Note that this operation preserves the properties of nowhere-zero and $k$-flow, so the existence of such flows depends only on the underlying signed graph and not the given orientation.

\begin{conjecture}[Bouchet \cite{Bouchet}]
If $G$ is a flow-admissable signed graph, then it has a nowhere-zero 6-flow.
\end{conjecture}

\begin{theorem}\label{thm:main}
If $G$ is a signed cyclically 5-edge-connected cubic graph that is flow-admissible, then $G$ has a nowhere-zero $8$-flow.
\end{theorem}

\begin{theorem}[Lu, Luo, and Zhang, \cite{no2edgedisjcycles}]\label{thm:no2disjNeg}
If $G$ is a flow-admissible signed graph with no two edge-disjoint negative cycles, then $G$ has a nowhere-zero $6$-flow.
\end{theorem}

When working with flows on ordinary (not signed) graphs, a classical result due to Tutte shows that a graph has a nowhere-zero $k$-flow if and only if it has a nowhere-zero flow in ${\mathbb Z}_k$, or in fact, more generally in any abelian group of order $k$.  One of the primary challenges in working with signed graphs is that this result no longer holds true.  There is a related approximate version due to Bouchet: If $\phi : E(G) \rightarrow \Z$ is an integer flow in the bidirected graph $G$, then there exists a $2k$-flow $\psi : E(G) \rightarrow \Z$ with $\psi \equiv \phi (\mod k)$.  This result is key in the proof of the 11 and 12 flow theorems, but this factor of 2 is more than we can afford in constructing our desired flow.  In our framework, we will be interested in converting a $\Z_3$-flow to an integer flow.  There is one case when this task is handled very nicely, namely the following theorem.  We give the proof because it is quick and instructive.

\begin{theorem}
\cite[Lemma 3.2]{XZ}
\label{Z3-to-3}
Let $G$ be a 2-connected bidirected cubic graph and let $\phi : E(G) \rightarrow {\mathbb Z}_3$ be a nowhere-zero flow.  Then there exists a 3-flow $\psi : E(G) \rightarrow {\mathbb Z}$ satisfying $\psi + 3 {\mathbb Z} = \phi$.
\end{theorem}

\begin{proof}
If $\phi(e) = -1$ we may adjust $e$ by reversing its orientation and setting $\phi(e) =1$; by this operation we may assume $\phi$ is constantly 1. This means that in our orientation, every vertex is either a source or sink.  Since $G$ is a 2-connected cubic graph, it contains a perfect matching $M$ and we may define $\psi : E(G) \rightarrow \Z$ by $\psi(e) = -2$ if $e \in M$ and $\psi(e) = 1$ otherwise.  It is immediate that $\psi$ is a nowhere-zero $3$-flow congruent to $\phi$.
\end{proof}

The assumption that $\phi$ is nowhere-zero is too strong for our purposes.  Instead our method relies upon the following closely related theorem which has a weaker conclusion (4-flow instead of 3-flow), but no longer requires the flow to be nowhere-zero.

\begin{theorem}\label{thm:pm}
Let $G$ be a 2-connected bidirected cubic graph and let $\phi : E(G) \rightarrow {\mathbb Z}_3$ be a flow\katso{with the property that $\{ e \in E(G) \mid \phi(e) = 0 \}$ is a matching}.  Then there exists a 4-flow $\psi : E(G) \rightarrow {\mathbb Z}$ satisfying $\psi + 3 {\mathbb Z} = \phi$.
\end{theorem}

This result is key to our approach, and the curious reader may skip ahead to Section~\ref{sec:4pre} to see the proof. As in the proof of Theorem \ref{Z3-to-3}, it relies upon the existence of a perfect matching, but in a certain alternate auxiliary graph.  For technical reasons to be described in the following section, we actually require a stronger form of the above theorem, and for this we need to prove the existence of perfect matchings in certain subcubic (i.e. maximum degree at most 3) graphs.  The following theorem gives a natural condition implying  the existence of a perfect matching that may be of independent interest, so we have included it in this introduction. \kat{It follows from Lemma \ref{perfect_matching}}.

\begin{theorem}
Let $G$ be obtained from a cyclically 5-edge-connected cubic graph by subdividing edges to form (degree 2) vertices $x_1, \ldots, x_{2k}$.  For every $1 \le i \le 2k$ let $C_i$ be a cycle of $G$ of odd length with $x_i \in V(C_i)$.  If $V(C_i) \cap V(C_j) = \emptyset$ whenever $i \neq j$ then $G$ has a perfect matching.
\end{theorem}

\section{Outline}

The proof of our main result is modeled after Seymour's original proof of his $6$-flow theorem \cite{Seymour6flow}, which we quickly recall. First, he reduces the problem to $3$-connected cubic graphs.  Define a \emph{generalized cycle} to be a graph that is either a cycle or a single isolated vertex.  For a 3-connected cubic graph $G$, he builds a list of generalized cycles, $C_1, \ldots, C_t$ so that $\{ V(C_1), \ldots, V(C_t) \}$ is a partition of $V(G)$ and so that for every $2 \le k \le t$ there are at least two edges between $V(C_k)$ and $\cup_{i=1}^{k-1} V(C_i)$.  He uses this structure to construct a $\Z_3$-flow $\phi : E(G) \rightarrow \Z_3$ so that $\phi(e) = 0$ only if $e \in \bigcup_{i=1}^t E(C_i)$.  By Tutte's Theorem, there exists a 3-flow $\psi : E(G) \rightarrow {\mathbb Z}$ congruent to $\phi$ modulo 3.  There is a 2-flow $\tau$ with support $\bigcup_{i=1}^t E(C_i)$ and then $2 \psi + \tau$ is a nowhere-zero 6-flow.  

Similar to Seymour, our proof of Theorem \ref{thm:main} will find a list of generalized cycles $C_1, \ldots, C_t$ in our signed graph $G$, together with a suitable function $\phi$ mapping the edges to $\Z_3$. However, the framework of signed graphs leads to numerous complications.  The first issue is that some of the $C_i$ cycles will be negative, and a negative cycle is not the support of a $2$-flow.  To deal with this, we will not be building our flow as a sum of two different flows, but instead we will build our flow as a sum of two functions that are close to being flows but have nonzero boundary at a small set of vertices.  Of course, we will need to arrange that these functions have opposite boundaries so that they sum to a flow.  \kat{To do so,}\katso{As a bookkeeping device,} we will modify our original graph by subdividing edges forming degree 2 vertices, and these new degree 2 vertices will be precisely those on which the boundaries of our functions are nonzero.  Let us introduce a key definition in accordance with this:  For a subcubic bidirected graph $G = (V,E)$ a function $\phi : E \rightarrow \Gamma$ is a \emph{preflow} if $\partial \phi (v) \neq 0$ if and only if $\mathrm{deg}(v) \in \{1,2\}$.  The final step in our argument will operate on a subdivision $G^*$ of the original graph where we combine a 2-preflow $\tau$ and a 4-preflow $\psi$ to get an 8-flow $\tau + 2 \psi$. 

We will then convert this $\Z_3$-preflow to a 4-preflow.  As mentioned in the introduction, this will require us to prove the existence of a perfect matching in a certain modified version of our graph.   The presence of this perfect matching would be immediate from Petersen's Theorem were our modified graph to be cubic, but the subdivisions we perform bring us instead to a subcubic graph.  Accordingly, our process for selecting generalized cycles is biased toward selecting positive cycles.  If we have selected the cycles $C_1, \ldots, C_k$ so far, then generally our next step will be to find a positive cycle $C_{k+1}$ so that at least two vertices in $V(C_{k+1})$ have a neighbour on a previously selected cycle.  When no such choice is possible (and our list does not yet cover the vertices) we prove that we may select a negative cycle as $C_{k+1}$ with the property that all but at most one vertex in $V(C_{k+1})$ has a neighbour on a previously selected cycle.  The structure of this selection process will be called upon to prove that the desired perfect matching exists.

There is another more immediate issue caused by the degree 2 vertices that we create by subdivision.  In order for the perfect matching we require to exist, it is necessary for the number of degree 2 vertices we create to be even. \kat{This boils down to forcing the number of even-length negative cycles to be even. When the remaining graph has no unbalanced theta, every negative cycle is its own block, and will eventually be a cycle in our list. When the remaining graph has an unbalanced theta, we are able to choose the parity of the number of negative cycles in our list. Thus we take special care when destroying the last unbalanced theta in the remaining graph. And this allows us to ensure the number of degree 2 vertices that we create is even.}

\kat{The remainder of the paper proceeds as follows. In Sections 3 and 4, we put together the tools to build our list of cycles, and finally develop a lemma to build this list in Section 5. In Section 6, we show that given the list of cycles, a suitable $\Z_3$-preflow exists. Then, using a perfect-matching lemma in Section 7, we show how to convert this $\Z_3$-preflow to a 4-flow in Section 8. Finally, in Section 9 we put all of the tools together to prove our main result.}

\section{Finding a usable cycle}

At a generic step of our cycle selection process we will have chosen a list of generalized cycles $C_1, \ldots, C_k$ and we will be interested in extending this list by choosing the next cycle $C_{k+1}$ from the graph $G' = G - \bigcup_{i=1}^k V(C_i)$.  If $G'$ contains a vertex $v$ with $\mathrm{deg}_{G'}(v) \le 1$ then this one-vertex subgraph will be a suitable generalized cycle.  Otherwise $G'$ will have vertices of degree 2 and 3 and we will attempt to extend our list by choosing a positive cycle $C_{k+1} \subseteq G'$ with the property that $C_{k+1}$ contains at least two vertices that have degree 2 in $G'$.  In this section we establish the basic lemma we will lean on to find such a cycle.  

We say that a generalized cycle $G$ in a subcubic signed graph is \emph{good} if $C$ is either a vertex of degree 0 or 1, or $C$ is a positive cycle containing at least two vertices of degree 2.  A generalized cycle is \emph{usable} if it is either good, or is a negative cycle containing at most one vertex not of degree 2.  Our goal in this section is to prove a lemma implying the existence of a usable cycle in any proper induced subgraph of a cyclically 5-edge-connected cubic graph.  We begin toward this with a well-known property of 3-connected graphs (and sketch the proof for completeness).  

\begin{lemma}
Let $G=(V,E)$ be a simple 3-connected graph with $|V| \ge 5$ and let $e \in E$.  Then either $G / e$ is 3-connected, or $G \setminus e$ is a subdivision of a 3-connected graph.
\end{lemma}

\begin{proof}[Proof sketch]
Let $e = xy$.  If $G/e$ is not 3-connected, it has a vertex cut of the form $\{x,y,z\}$.  If $G \setminus e$ is not a subdivision of a 3-connected graph, it has a vertex cut, say $\{u,v\}$, splitting $V \setminus \{u,v\}$ into two sets $X,Y$ with $|X|, |Y| \ge 2$ so that $x \in X$ and $y \in Y$.  Together these properties contradict the assumption that $G$ is 3-connected.
\end{proof}

Our next lemma calls on the previous to prove the existence of theta subgraphs that will be of use. \kat{Let $V_i(G)$ be the set of vertices of degree $i$ in $G$.}

\begin{lemma}\label{lem:theta}
Let $G$ be obtained from a 3-connected graph by doing $k$ subdivisions.  Assume that $k \ge 4$ and that no edge of the original graph is subdivided more than $k-2$ times.  Then $G$ either contains a theta subgraph so that each branch contains a vertex in $V_2(G)$ or $G$ has a theta subgraph so that two branches have at least two vertices in $V_2(G)$.  
\end{lemma}

\begin{proof}
We proceed by induction on $|V| + |E|$.  If $G$ contains $x,y \in V_2(G)$ with $x,y$ adjacent, then by our assumptions we may choose a cycle $C \subseteq G \setminus x$ containing at least two vertices in $V_2(G)$.  Now $C$ can be completed to a theta containing $x$, and this theta satisfies the lemma.  Accordingly, we may now assume that $V_2(G)$ is an independent set.  The result is straightforward to verify if either $G$ is a subdivision of $K_4$ or if $V_{3+}(G)$ is an independent set, so we may assume neither of these hold.  Choose an edge $e$ so that both ends of $e$ are in $V_{3+}(G)$.  By the previous lemma, either $G / e$ or $G \setminus e$ is a subdivision of a 3-connected graph, and in both cases the result follows by applying induction to this graph (suppressing any newly created degree two vertices).
\end{proof}

With this we are ready to prove our main lemma from this section.  

\begin{lemma}
\label{usable_cycle}
Let $G = (V,E)$ be obtained from a cyclically 5-edge-connected cubic graph by deleting at least three vertices.  Then $G$ contains a usable generalized cycle.
\end{lemma}

\begin{proof}
We may assume that $G$ is connected, and not a generalized cycle (in which case $G$ itself is usable).  Since $V_2(G) \neq \emptyset$, we may choose a minimal proper subset $X \subseteq V(G)$ so that $G[X]$ contains a cycle and $e(X,\overline{X}) \le 2$.  If $e(X,\overline{X}) = 1$ and $G[X]$ is a cycle, then this is a usable cycle and we have nothing left to prove.  Otherwise in this case we let $H$ denote the graph obtained from $G[X]$ by suppressing the unique vertex of $X$ that had a neighbour in $\overline{X}$.  If $e(X,\overline{X}) = 2$ then it follows from our connectivity assumptions on $G$ that there exists $y \in V_2(G) \setminus X$ and we may choose a  minimal path $P$ containing $y$ so that both ends of $P$ are in $X$.  Let $H$ be obtained from $P \cup G[X]$ by suppressing all interior vertices of $P$ apart from $y$.  It follows from Lemma \ref{lem:theta} that $H$ contains a theta with at least one degree 2 vertex on each branch, or at least two degree 2 vertices on at least two branches.  The original graph $G$ must also contain such a theta, and since every theta contains a positive cycle, this implies the existence of a good cycle, as desired.
\end{proof}

\section{Using an unbalanced theta}\label{unbaltheta}

As described earlier, we will eventually need to find a perfect matching in a certain (modification of) a subdivision of our original graph $G$.  This condition in turn forces a certain parity property on our list of cycles.  In this section we develop the main tool to ensure that this condition is satisfied.  Before proving this lemma, we introduce a basic observation that we will call upon in its proof.

\begin{observation}
Let $G$ be a signed graph, let $H \subseteq G$ be connected, and assume that both $H$ and $G - E(H)$ are balanced.  If $G$ is not balanced, there exist paths $P \subseteq H$ and $Q \subseteq G-E(H)$ so that $P \cup Q$ is a negative cycle.
\end{observation}

\begin{proof}
By possibly modifying the signature, we may assume that every edge in $E(G) \setminus E(H)$ is positive.  Since $H$ is balanced, we may express $V(H)$ as a disjoint union $U_0 \cup U_1$ so that all edges of $H$ with one end in $U_0$ and the other in $U_1$ are negative and all other edges of $H$ are positive.  If there exists a path in $G - E(H)$ from $U_0$ to $U_1$, then observation follows by choosing a minimal such path $Q$ together with a path $P \subseteq H$ with the same ends.  If no such path exists, then $\delta_H(U_0)$ is an edge-cut of $G$ and switching on this edge-cut gives an equivalent signature of $G$ with all edges positive, thus showing that $G$ is balanced.
\end{proof}

\begin{lemma}
\label{fragile}
Let $G = (V,E)$ be a 2-connected graph obtained from a cyclically 5-edge-connected cubic graph by deleting a set of at least three vertices.  Assume that $G$ contains an unbalanced theta, but that $G - E(C)$ does not contain an unbalanced theta for every good cycle $C$.  Then there exist two cycles $C_0,C_1 \subseteq G$ satisfying the following properties:
\begin{enumerate}
\item $C_1$ is negative and $C_0$ is good (and thus positive)
\item $G - E(C_i)$ does not have an unbalanced theta for $i=0,1$
\item $G - E(C_0)$ and $G - E(C_1)$ have the same negative cycles.  
\end{enumerate}
\end{lemma}

\begin{proof}
We split the argument into cases.

\medskip

\noindent{Case 1:} Every good cycle intersects every negative cycle.  

\medskip

Choose a good cycle $C$, choose a vertex $x \in V(C)$ with $\mathrm{deg}_G(x) = 2$ and apply the observation to choose a path $P \subseteq G \setminus E(C)$ so that $P \cup C$ is an unbalanced theta.  Let $Q$ denote the path in $C$ with the same ends as $P$ not containing $x$, and assume that among all possible choices, $P$ has been selected so as to minimize $Q$.  Let $C'$ denote the cycle in $C \cup P$ not containing the path $Q$ and note that $C'$ is negative.  We claim that $G - E(C')$ is balanced.  To see this, suppose for a contradiction that $N \subseteq G - E(C')$ is a negative cycle.  Note that $N$ must contain an edge in $C$ since $G - E(C)$ is balanced.  Now applying the observation to the graph $C \cup N$ and the subgraph $H = C$ gives us a path $P' \subseteq N$ so that $C \cup P'$ is an unbalanced theta, and this path $P'$ contradicts our choice of $P$.  Thus $G - E(C')$ is balanced, as desired, and now the cycles $C$ and $C'$ satisfy the lemma.

\medskip

\noindent{Case 2: } There does not exist an induced negative cycle $N \subseteq G$ so that $d(V(N)) = 2$.  

\medskip

We shall assume the first case is not satisfied, and choose a good cycle $C$ and a negative cycle $N$ with $V(C) \cap V(N) = \emptyset$.  Choose two degree two vertices $u,v \in V(C)$ and let $Q_1,Q_2$ be the two paths in $C$ from $u$ to $v$.  Let $H$ be the component of $G - E(C)$ containing $N$ and note that $N$ is a block of $H$ (otherwise $H$ would have an unbalanced theta).  By the assumption for this case, there exist (at least) three vertices $x_1, x_2, x_3 \in V(N)$ that have a neighbour not in $V(N)$ and for $i=1,2,3$ we may choose a path $P_i \subseteq H \setminus E(N)$ from $x_i$ to $V(C)$.  Assume that $P_i$ has ends $x_i,y_i$ for $1 \le i \le 3$ and that $P_i$ is internally disjoint from $V(C)$.  We may assume (without loss) that $P_1$ and $P_2$ both end in $V(Q_1)$ and we let $R$ denote the subpath of $Q_1$ from $y_1$ to $y_2$.  Let us assume that our $P_i$ paths have been chosen so that $R$ is minimal, and note that this implies that $R$ is internally disjoint from $V(H)$.  Now the graph $F$ obtained from $C \cup P_1 \cup P_2 \cup N$ by removing the edges and internal vertices of $R$ is an unbalanced theta.  It has two cycles, say $C_1$ and $C_2$, containing $u$ and $v$ with one positive (and thus good) and one negative.  Both $G- E(C_1)$ and $G-E(C_2)$ will have components that are subgraphs of $H$, but apart from these, every component of $G-E(C_i)$ is also a component of  $G-E(C_{3-i})$.  It follows that $C_1$ and $C_2$ satisfy the lemma.

\medskip

\noindent{Case 3: } There exists an induced negative cycle $N \subseteq G$ so that $d(V(N)) = 2$.

\medskip

Let $x,y$ be the two vertices in $V(N)$ that have a neighbour in $V \setminus V(N)$, let $P_1, P_2 \subseteq N$ be the two paths in $N$ from $x$ to $y$, and let $G' = G - E(N)$.  Since $G$ is a subgraph of a cyclically 5-edge-connected cubic graph we must have $|V(N)| \ge 5$ and thus we may assume, without loss what $|V(P_1)| \ge 4$ (so $P_1$ contains at least two vertices of degree 2).  If there is a path $Q \subseteq G-E(N)$ from $x$ to $y$ so that $P_1 \cup Q$ is positive, then $P_1 \cup Q$ is good and this cycle together with $P_2 \cup Q$ satisfy the lemma.  Accordingly, we may assume no such path exists.  It then follows that $G-E(N)$ is balanced (otherwise, a negative cycle, say $N'$, together with two vertex disjoint paths from $N'$ to $\{x,y\}$ would give a suitable path $Q$).  Therefore, by possibly switching to an equivalent signature, we may assume that there is just one negative edge in $G$ and it lies on the path $P_1$.  

To complete the proof, it suffices to find a good cycle $C$ containing $x$; in this case $G \setminus V(C)$ and $G \setminus N$ are both balanced and together satisfy the lemma.  If $G' = G \setminus V(N)$ contains a cut-edge $e$, then choose $u,v \in V_2(G)$ in separate components of $G' \setminus e$.  There exists a path $P \subseteq G$ from $x$ to $y$ containing both $u,v$ and now $P_2 \cup P$ is a good cycle containing $x$.  Accordingly, we may assume $G'$ is 2-connected.  Choose $u,v,w \in V_2(G) \cap V(G')$, let $H$ be a minimally 2-connected subgraph of $G'$ containing $u,v,w$ and note that $H$ is either a cycle or a theta subgraph.  Choose two vertex disjoint paths $Q_1, Q_2$ in $G$ from $V(N)$ to $V(H)$.  In all cases (because there is no good cycle disjoint from an unbalanced theta), $P_2 \cup Q_1 \cup Q_2 \cup H$ contains a good cycle containing $x$, thus completing the proof.
\end{proof}

\section{Constructing a sequence of generalized cycles}

In this section, we show how to find the desired sequence of generalized cycles.  We require one preliminary lemma for this allowing us to choose the first cycle for the sequence.

\begin{lemma}
Let $G$ be a 3-connected signed graph that has two vertex disjoint negative cycles.  Then one of the following holds:
\begin{enumerate}
\item $G$ has vertex disjoint negative cycles $C_1$, $C_2$ so that $V(C_1) \cup V(C_2) =V(G)$, or
\item $G$ has a negative cycle $C$ so that $G \setminus V(C)$ has an unbalanced theta.
\end{enumerate}
\end{lemma}

\begin{proof}
Choose two vertex disjoint negative cycles, say $C_1$ and $C_2$.  If the first condition is satisfied, then we are finished.  Otherwise, there is a vertex $v \in V(G) \setminus (V(C_1) \cup V(C_2))$.  By Menger's theorem we may choose three paths $P_1, P_2, P_3$ all starting at $v$ and ending in $V(C_1) \cup V(C_2)$ so that these paths are disjoint apart from $v$.  We may assume further that each $P_i$ is internally disjoint from $V(C_1) \cup V(C_2)$.  We may assume without loss that $P_1$ and $P_2$ end in $V(C_2)$.  Now $C_1$ is a cycle that satisfies the second condition as $P_1 \cup P_2 \cup C_2$ is an unbalanced theta.
\end{proof}

We are now ready for the main lemma.

\begin{lemma}
\label{cycle_list}
Let $G$ be a signed cyclically 5-edge-connected cubic graph with two disjoint negative cycles.  Then there exists a sequence of generalized cycles $C_1, \ldots, C_t$ and a division of the negative cycles in this list into \emph{special} and \emph{ordinary} so that the following properties are satisfied:
\begin{enumerate}
\item $C_1, \ldots, C_t$ are vertex disjoint and $V(G) = \bigcup_{i=1}^t V(C_i)$,
\item $C_1$ is special, 
\item At most one cycle other than $C_1$ is special, 
\item Every ordinary cycle is induced,
\item Every positive $C_k$ satisfies $e( V(C_k), \cup_{i=1}^{k-1} V(C_i) ) \ge 2$
\item Every ordinary $C_k$ satisfies $e( V(C_k), \cup_{i=k+1}^{t} V(C_i) ) \le 1$
\item The number of even length ordinary cycles plus the number of special cycles is even.
\end{enumerate}
\end{lemma}

\begin{proof}
We begin by applying the first lemma from this section.  If we get the first outcome, then these two cycles are declared special and then they satisfy the lemma, completing the proof in this case.  If we get the second outcome, then we take this cycle to be $C_1$, declare it special, and this will be the start of our sequence.  In the first part of our selection process, we will grow our list of cycles while maintaining the property that the graph obtained by deleting the vertices of all chosen cycles still contains an unbalanced theta.  Suppose we have chosen $C_1, \ldots, C_k$ and $G \setminus (\bigcup_{i=1}^k V(C_i) )$ contains an unbalanced theta.  It follows from Lemma \ref{usable_cycle} that we may choose a usable cycle $C_{k+1}$ from $G \setminus (\bigcup_{i=1}^k V(C_i) )$.  If we can choose such a cycle $C_{k+1}$ so that $G \setminus (\bigcup_{i=1}^{k+1} V(C_i) )$ contains an unbalanced theta, then we add this cycle to the list, declaring it to be ordinary if $C_{k+1}$ is negative, and we carry on.  Note that if $C_{k+1}$ is positive, then it will satisfy condition 5 (since it was a good cycle), and if $C_{k+1}$ is ordinary, then it will satisfy conditions 4 and 6.  If no such cycle exists, then the cycle $C_{k+1}$ will be chosen so that it is either a good cycle, or a negative cycle that is declared to be special.  We may apply Lemma \ref{fragile} to choose such a cycle $C_{k+1}$ so that $G \setminus (\bigcup_{i=1}^{k+1} V(C_i) )$ has no unbalanced theta and the following condition is satisfied:
\begin{itemize}
\item[$(\star)$] The number of special cycles in $C_1, \ldots, C_{k+1}$ plus the number of ordinary cycles of even length in $C_1, \ldots, C_{k+1}$ plus the number of negative cycles of even length in $G \setminus (\bigcup_{i=1}^{k+1} V(C_i) )$ is even.
\end{itemize}
Having done this, we will carry on applying Lemma \ref{usable_cycle} to choose usable cycles from the remaining graph (declaring any negative cycle selected to be ordinary), completing the list $C_1, \ldots, C_t$ so that $V(G) = \bigcup_{i=1}^t V(C_i)$.  Again note that the definition of usable implies that any positive cycle selected will satisfy condition 5 and any ordinary cycle selected will satisfy 4 and 6.  The key observation we need to complete the proof is that any negative cycle $N$ in the graph $G \setminus (\bigcup_{i=1}^{k+1} V(C_i) )$ is a block of this graph and must therefore be selected as an ordinary cycle $C_i$ at some point in this process.  It follows from this observation and $(\star)$ that our resulting list of cycles $C_1, \ldots, C_t$ satisfies condition 7.  All of the other conditions are satisfied by our process, and this completes the proof.
\end{proof}

\section{Forming a ${\mathbb Z}_3$-preflow}

In this section we will use the list of cycles developed in the previous section to construct a suitable ${\mathbb Z}_3$-preflow.  We begin with two straightforward lemmas concerning basic properties of signed cycles.  

\begin{lemma}
\label{cycle_boundary1}
Let $C$ be a bidirected cycle and let $b: V(C) \rightarrow {\mathbb Z}_3$.  If one of the following holds:
\begin{itemize}
\item $C$ is negative, or
\item All edges of $C$ are positive and $\sum_{v \in V(C)} b(v) = 0$,
\end{itemize}
Then there exists $\tau : E(C) \rightarrow {\mathbb Z}_3$ so that $\partial \tau = b$.
\end{lemma}

\begin{proof}
For the first part, observe that for every $u \in V(C)$ there is a function $\tau_u : E(C) \rightarrow \{ -1,1\}$ so that $\partial \tau_u(v) = 0$ for every $v \in V(C) \setminus \{u\}$ and $\partial \tau_u(u) \neq 0$.  A suitable combination of $\tau_u$ functions will give a function $\tau$ with any desired boundary.  For the second part, choose $u \in V(C)$ and note that we may greedily choose a function $\tau : E(C) \rightarrow {\mathbb Z}_3$ so that $\partial \tau(v) = b(v)$ holds for every $v \in V(C) \setminus \{u\}$.  Since every edge of $C$ is positive we have $\sum_{v \in V(C)} \partial \tau(v) = 0 = \sum_{v \in V(C)} b(v)$ and it follows that $\partial \tau = b$ as desired.
\end{proof}

\begin{lemma}
\label{cycle_boundary2}
If $C$ is a bidirected cycle that is negative and $|V(C)|$ is odd, then there exists $\tau : E(C) \rightarrow \{-1,1\}$ so that $\partial \tau(v) \neq 0$ for every $v \in V(C)$.
\end{lemma}

\begin{proof}
By possibly reorienting $C$ we may assume every vertex of $C$ is either a source or a sink.  Now we may take $\tau$ to be the constant 1.
\end{proof}

If $G$ is a signed graph, and $e=uv \in E(G)$, then to \emph{subdivide} $e$ is to replace $e$ in $G$ by a path $P$ of length 2, with the same ends as $e$, so that $\sigma(P) = \sigma(e)$. We say a signed graph $G'$ is a subdivision of $G$ if $G'$ can be obtained from $G$ by a sequence of subdivisions.

\begin{lemma}
\label{construct_preflow}
Let $G$ be a bidirected cubic graph and let $C_1, \ldots, C_t$ be a list of cycles in~$G$ satisfying the conditions from Lemma \ref{cycle_list}.  Let $G^*$ be a bidirected graph obtained from $G$ by subdividing one edge from each cycle in this list that is either special or is ordinary with even length, and for every $C_i \subseteq G$ let $C_i^*$ be the corresponding cycle in $G^*$.  Then there exists a preflow $\phi : E(G^*) \rightarrow {\mathbb Z}_3$ satisfying:
\begin{itemize}
\item Every edge $e$ with $\phi(e) = 0$ lies on a cycle $C_i^*$ where $C_i$ is either positive or special.
\end{itemize}
\end{lemma}

\begin{proof}
For $1 \le k \le t$ let $U_k =  \bigcup_{i=k}^t V(C_i^*)$.  The cycles $C_1, \ldots, C_t$ are classified as positive, special, or ordinary, and we use the corresponding classification for $C_1^*, \ldots, C_t^*$.  We may assume by possibly changing to an equivalent signature that every positive $C_i^*$ has all edges positive.  We begin by setting $\phi : E(G^*) \rightarrow {\mathbb Z}_3$ to be the constant 0.  We will now adjust $\phi$ in steps starting from the back of the list of cycles.  Maintaining the property that after we have adjusted $\phi$ on all edges incident with $C_k^*$, then $\partial \phi(v)$ is nonzero if and only if $v \in V_2(G)$.  If we have already adjusted $\phi$ on the edges incident with a vertex in $U_{k+1}$, then we next adjust $\phi$ on edges incident with a vertex in $V(C_{k}^*)$ but not incident with a vertex in $U_{k+1}$ according to the following process: 

\medskip

\noindent{Case 1: } $C^*_{k}$ is positive

\medskip

Define $\phi(e) = 1$ for every chord edge $e$ of $C^*_{k}$.  Define $F = E( V(C^*_k), \overline{U_k} )$ and note that by assumption $|F| \ge 2$.  It follows that we may assign values $\phi(f) = \pm 1$ to every $f \in F$ so that $\sum_{v \in V(C^*_k)} \partial \phi(v) = 0$.  Lemma \ref{cycle_boundary1} implies that we may choose $\tau : E(C^*_k) \rightarrow {\mathbb Z}_3$ so that $\partial \tau(v) = - \partial \phi(v)$ holds for every $v \in V(C^*_k)$.  Now setting $\phi = \phi + \tau$ ensures that $\phi$ has zero boundary at every vertex in $V(C_k^*)$, as desired.

\medskip

\noindent{Case 2: } $C^*_{k}$ is ordinary

\medskip

Note that thanks to our subdivisions, the cycle $C^*_{k}$ has odd size.  Apply Lemma \ref{cycle_boundary2} to choose $\tau : V(C^*_k) \rightarrow {\mathbb Z}_3$ so that $\partial \tau (v)$ is nonzero for every $v \in V(C^*_k)$.  Since $C^*_{k}$ is ordinary, there is at most one vertex $v \in V(C^*_k)$ adjacent to a vertex in $U_{k+1}$; if such a vertex exists, we call it $w$.  By possibly replacing $\tau$ by its inverse, we may assume $\partial \tau(w) = - \partial \phi(w)$ if $w$ exists.  Now set $\phi = \phi + \tau$ and extend the domain of $\phi$ to include $E(V(C^*_k), \overline{U_k} )$ by choosing $\phi(e) = \pm 1$ for every such edge so that the resulting function $\phi$ has zero boundary at every degree 3 vertex in $V(C_k^*)$.  

\medskip

\noindent{Case 3: } $C^*_k$ is special.

\medskip

Define $\phi(e) = 1$ for every chord edge $e$ of $C^*_{k}$ and for every edge in $E( V(C^*_k), \overline{U_{k}} )$.  It follows from Lemma \ref{cycle_boundary1} that we may choose $\tau : E(C^*_k) \rightarrow {\mathbb Z}_3$ so that $\partial \tau(v) = - \partial \phi(v)$ holds for every $v \in V(C^*_k)$ of degree 3 and $\partial \tau(w) = 1$ for the unique vertex $w \in V(C^*)$ of degree 2.  Now we set $\phi = \phi + \tau$.

\medskip

It follows from this process that our final function $\phi$ is a ${\mathbb Z}_3$-preflow with the property that every edge $e$ with $\phi(e) = 0$ lies on a positive or special cycle, as desired.
\end{proof}

\section{Perfect matchings in subcubic graphs}\label{sec:pm}

In this section we will prove a lemma giving the existence of a perfect matching in a subcubic graph under certain assumptions.  Our result relies upon the following famous result.  For a graph $G$ we let $\mathrm{odd}(G)$ denote the number of components of $G$ with an odd number of vertices.

\begin{theorem}[Tutte \cite{tuttepm}]
A graph $G$ has a perfect matching if and only if there does not exist a set $X \subseteq V(G)$ so that $|X| < \mathrm{odd}(G \setminus X)$.
\end{theorem}

Let $G = (V,E)$ be a connected graph and let $X \subseteq V$.  We say that a cycle $C$ \emph{straddles} the edge-cut $\delta(X)$ is $C$ contains an edge with both ends in $X$ and an edge with both ends in $V \setminus X$.

\begin{lemma}
\label{perfect_matching}
Let $G$ be a subdivision of a 3-connected cubic graph, let $x_1, \ldots, x_k \in V_2(G)$, and for $1 \le i \le k$ let $C_i$ be a cycle of odd length in $G$ with $x_i \in V(C_i)$.  If the following conditions are satisfied:
\begin{itemize}
\item $|V_2(G)|$ is even,
\item $V(C_i) \cap V(C_j) = \emptyset$ whenever $i \neq j$,
\item $|V_2(G) \setminus \{x_1, \ldots, x_k\}| < 6$,
\item No 3-edge-cut $S$ is straddled by $C_i$ for some $1 \le i \le k$, and
\item No 4-edge-cut $S$ is straddled by $C_i$ and $C_j$ for some $1 \le i < j \le k$,
\end{itemize}
then $G$ contains a perfect matching.
\end{lemma}

\begin{proof}
Suppose (for a contradiction) that $G$ does not have a perfect matching and apply Tutte's Theorem to choose a maximal set $X \subseteq V$ so that $|X| < \mathrm{odd}(G \setminus X)$.  It follows from the maximality of $X$ that every component of $G \setminus X$ has odd size (otherwise move a vertex from an even-size component of $G \setminus X$ to the set $X$ to enlarge it).  Note that since our assumptions imply $|V(G)|$ is even, $|X|$ and $\mathrm{odd}(G \setminus X)$ have the same parity, and thus  $\mathrm{odd}(G \setminus X) \ge |X| + 2$.  Let ${\mathcal J} = \{1 \le i \le k \mid \mbox{$C_i$ has an edge with both ends in $X$} \}$.   Then we have the following upper bound:
\begin{align*}
E(X,\overline{X})	&=		\sum_{x \in X} \mathrm{deg}(x) - 2 e(X)  	\\
				&\le 		3 |X| - |{\mathcal J}|. 
\end{align*}
Now consider a component $H$ of $G \setminus X$.  We say that $H$ is \emph{negligible} if $V(H) \subseteq V_2(G) \setminus \{x_1, \ldots, x_k\}$.  We let ${\mathcal N}$ denote the set of negligible components and note that by our assumptions $| {\mathcal N}| \le 5$, and each such component $H$ will have $e(V(H),X) \ge 2$.  We say that $H$ is \emph{critical} if $H$ consists of a single vertex that is one of $x_1, \ldots, x_k$ (in which case $e(V(H),X) = 2)$ and we let ${\mathcal C}$ denote the set of critical components.  Finally, we call $H$ \emph{basic} if it is not negligible or critical, and we let ${\mathcal B}$ denote the set of basic components.  Every basic component $H$ has $e(V(H),X) \ge 3$, but we get a better bound by considering straddling cycles.  It follows from our assumptions that every $H \in {\mathcal S}$ has the property that $e( V(H), X)$ is at least 3 plus the number of cycles $C_i$ that straddle the cut $\delta_G( V(H) )$.  Let ${\mathcal S}$ denote the set of all those cycles $C_i$ for which $C_i$ straddles an edge cut of the form $\delta_G(V(H))$ where $H \in {\mathcal B}$.  Now we have the lower bound:
\begin{align*}
E(X,\overline{X})	
	&=		\sum_{H \in \mathrm{odd}(G \setminus X) }	e( V(H), X)	\\
	&\ge		3 | \mathrm{odd}(G \setminus X) | - | {\mathcal N} | - | {\mathcal C} | + | {\mathcal S} |	\\
	&\ge		3 |X| + 1 - | {\mathcal C} | + | {\mathcal S} |
\end{align*}
Combining our bounds on $E(X,\overline{X})$ gives us the inequality 
\[ | {\mathcal C} |  \ge  1 + |{\mathcal S}| + | {\mathcal J} | \]
but this is a contradiction.  By parity, every time a vertex $x_i$ forms a critical component, the corresponding cycle $C_i$ has odd length and must therefore either contain an edge with both ends in $X$ or must straddle an edge cut of the form $\delta_G(V(H))$ where $H \in {\mathcal B}$.  This contradiction completes the proof.
\end{proof}

\section{$\Z_3$-preflows to 4-preflows}\label{sec:4pre}

As mentioned in the introduction, a general difficulty in the subject of flows in signed graphs is the fact that Tutte's theorem for ordinary graphs which converts from a modular flow to an integer flow is no longer true in general in the setting of signed graphs.  In this section we will show how to convert a $\Z_3$-flow in a 2-connected signed \kat{sub}cubic graph to a (congruent) 4-flow.  In fact, for reasons described in the previous section, we will need to operate in the setting of preflows on subcubic graphs.  

Let $\phi : E(G) \rightarrow \Z_3$ be a pre-flow in a bidirected graph $G$ with all vertices of degree 2 or 3, and assume that $Z = \{ e \in E(G) \mid \phi(e) = 0 \}$ is a matching. By possibly changing to an equivalent signature, we may assume that every edge in $Z$ is positive.  If $e \in E(G)$ has $\phi(e) = -1$ then we may reverse the orientation of $e$ and set $\phi(e) = 1$ while maintaining our preflow.  Doing this for all such edges gives a preflow $\phi$ with $\phi(E(G)) \subseteq \{0,1\}$.  Next modify $G$ to form a new bidirected graph $G^+$ by adding a single new vertex $w$ and then for every vertex $v \in V(G)$ with $\mathrm{deg}(v) = 2$ we add a positive edge with ends $v,w$.  We may extend the function $\phi$ to a function $\phi^+$ with domain $E(G^+)$ by assigning each newly added edge $e$ a suitable orientation and setting $\phi^+(e) = 1$ so that $\partial \phi^+(v) = 0$ holds for every $v \in V(G)$ (but possibly not at the vertex $w$).  

Let $e = uu' \in E(G^+)$ have $\phi(e) = 0$, let $h_1,h_2$ $(h_1',h_2')$ be the half edges incident with $u$ $(u')$ not contained in $e$, and let $f_1,f_2$ $(f_1',f_2')$ be the edges containing $h_1,h_2$ $(h_1',h_2')$ respectively.  It follows from our orientation that we may assume $h_1,h_1'$ are positive while $h_2,h_2'$ are negative.  Modify the graph by subdividing $f_1$ and $f_2'$ and adding a new edge in between.  Let $H$ be the graph (not signed or oriented) obtained by doing this process at every edge where $\phi$ is zero, and then deleting the vertex $w$.  We call $H$ the \emph{auxiliary} graph \emph{associated with} $G$ and $\phi$.  

\begin{lemma}\label{4preflow}
Let $G$ be a bidirected graph with all vertices of degree 2 or 3, let $\phi : E(G) \rightarrow {\mathbb Z}_3$ be a preflow with the property that $Z = \{ e \in E(G) \mid \phi(e) =0 \}$ is a matching.  If the auxiliary graph associated with $G$ and $\phi$ has a perfect matching, then there exists a 4-preflow $\psi : E(G) \rightarrow \Z$ so that
\begin{itemize}
\item $\psi + 3 \Z = \phi$
\item $\partial \psi (v) = \pm 1$ for every $v \in V_2(G)$.
\end{itemize}
\end{lemma}

\begin{proof}
Assume that the signature has been chosen so that all edges of $Z$ are positive and $\phi$ and the orientation have been adjusted so that $\phi(E(G)) \subseteq \{0,1\}$.  Let $H$ be the associated auxiliary graph and let $M$ be a perfect matching in $H$.  For every edge $e \in E(G)$ with $\phi(e) \neq 0$, if $e \in M$ then we set $\psi(e) = -2$ and if $e \not\in M$ we set $\psi(e) = 1$.  Next suppose that $e = uu' \in E(G)$ has $\phi(e) = 0$.  Let $f_1,f_2'$ be the edges incident with $u, u'$ that were subdivided in constructing $H$ and let $e'$ be the new edge added between the degree 2 vertices resulting from this subdivision (so $e,f_1,e',f_2'$ form the edges of a 4-cycle).  If both $e,e'$ or neither $e,e'$ are contained in $M$ then we set $\psi(e) = 0$. \kat{Otherwise exactly one of $e,e'$ is in $M$; set $\psi(e) = \pm3$, choosing the proper sign so that $\partial\psi(u) = 0 = \partial\psi(u')$.}
\end{proof}

\kat{Note that Lemma \ref{4preflow} implies Theorem \ref{thm:pm}. Let $G$ be a 2-connected bidirected cubic graph, and $\phi: E(G) \to \Z_3$ be a flow. If $Z = \{ e \in E(G) \mid \phi(e) =0 \}$ is not a matching, then by definition of flow there is a vertex $x \in V(G)$ incident to three edges in $Z$. Construct $G'$ from $G$ by replacing every such $x$ with a positive cycle $C$ of length $3$, so that every vertex in $C$ has degree $3$. And extend $\phi$ to $E(C)$ so that it remains a pre-flow and $\phi(e) = \pm1$ for every $e \in E(C)$. Then $Z$ is a matching in $G'$. Apply Lemma \ref{4preflow} to $G'$ and the extended $\phi$. Note that the auxiliary graph associated with $G'$ and $\phi$ is 2-connected and cubic, and so has a perfect matching by Petersen's Theorem. Thus there is a $4$-flow $\psi: E(G') \to \Z$ satisfying $\psi + 3\Z = \phi$. And $\psi$ is a flow when restricted to $G$ because the net flow across any edge-cut with one side balanced is zero.}

\section{Proof of the main theorem}

In this section we put together our lemmas to complete the proof of our main theorem.

\begin{proof}[Proof of Theorem \ref{thm:main}]
Let $G$ be a bidirected cyclically 5-edge-connected cubic graph.  If $G$ does not have two vertex disjoint negative cycles, then $G$ has a nowhere-zero 6-flow by Theorem \ref{thm:no2disjNeg} and we have nothing left to prove.  Otherwise, apply Lemma \ref{cycle_list} to choose a list of vertex disjoint cycles $C_1, \ldots, C_t$ and modify $G$ to form $G^*$ by subdividing one edge from each $C_i$ that is either special or is ordinary with even length. As before, we let $C_1^*, \ldots, C_t^*$ be the corresponding list of cycles in $G^*$.  Next apply Lemma \ref{construct_preflow} to choose a  preflow $\phi : E(G^*) \rightarrow {\mathbb Z}_3$.  By reorienting, we may assume that $-1$ is not in the range of $\phi$.  Note furthermore that by construction $\{ e \in E(G^*) \mid \phi(e) = 0 \}$ is a matching.  Let $H$ be the auxiliary graph associated with $G$ and $\phi$.  Consider an ordinary cycle $C_i$ of even length.  In forming the graph $G^*$ one edge of $C_i$ was subdivided one time with a new vertex, say $x_i$.  This means that the cycle $C_i^*$ will have odd length.  It follows from Lemma \ref{construct_preflow} that every edge $f$ incident with a vertex in $V(C_i^*)$ has $\phi(f) \neq 0$ and this means that these edges are not involved in any gadgets used in constructing our auxiliary graph $H$.  Therefore $C_i^*$ is an odd length cycle in $H$.  Since there are at most two degree two vertices that do not lie on ordinary cycles, it now follows from lemma \ref{perfect_matching} that $H$ contains a perfect matching.  Applying Lemma \ref{4preflow} gives us a 4-preflow $\psi$ of $G^*$ with the properties that $\psi(e) = 0$ only if $e$ lies on a positive or special cycle, and $\partial \psi(v) = \pm 1$ for every $v \in V_2(G^*)$.  Finally we construct a matching 2-preflow $\tau : E(G^*) \rightarrow \{-1,0,1\}$ starting from $\tau = 0$ as follows.  For every $C_i$ that is either special, or ordinary and even length, let $x_i$ be the degree 2 vertex created when subdividing an edge of $C_i$ to form $C_i^*$.  Modify $\tau$ by assigning $\tau(e) = \pm 1$ for every $e \in E(C_i^*)$ so that $\partial \tau(v) = 0$ for every $v \in V(C_i^*) \setminus \{x_i\}$ and so that $\partial \tau( x_i) = -2 \partial \psi(x_i)$.  It now follows that $\tau + 2 \psi$ is a nowhere-zero 8-flow of $G^*$, and this gives the desired flow in $G$, thus completing the proof.
\end{proof}

\bibliographystyle{acm}
\bibliography{bib}

\end{document}